\documentclass[12pt]{amsart}

\makeatletter
\@namedef{subjclassname@2020}{\textup{2020} Mathematics Subject Classification}
\makeatother

\newcommand{\tf}[2]{\tfrac{#1}{#2}}
\newcommand{\tff}[2]{\tfrac{#1}{#2}}

\newcommand{\Li}{\te{Li}}

\newcommand{\qte}[1]{\q\te{#1}}

\newcommand{\tri}{\triangle}

\newcommand{\dsum}{\di\sum}
\newcommand{\tsum}{\textstyle\sum}

\newcommand{\wi}{\ensuremath{W_{I}}}

\renewcommand{\wi}[1]{\di\int_{0}^{1}{#1}\,\ff{du}{\rt{1-u^{2}}}}

\DeclareFontFamily{U}{mathx}{\hyphenchar\font45}
\DeclareFontShape{U}{mathx}{m}{n}{
      <5> <6> <7> <8> <9> <10>
      <10.95> <12> <14.4> <17.28> <20.74> <24.88>
      mathx10
      }{}
\DeclareSymbolFont{mathx}{U}{mathx}{m}{n}
\DeclareFontSubstitution{U}{mathx}{m}{n}
\DeclareMathAccent{\widecheck}{0}{mathx}{"71}

\newcommand{\inrr}{\ensuremath{\in\rr}}

\renewcommand{\kill}[1]{}
\newcommand{\dummy}[1]{\mbox{}}

\makeatletter
\newcommand{\xequal}[2][]{\ext@arrow 0055{\equalfill@}{#1}{#2}}
\def\equalfill@{\arrowfill@\Relbar\Relbar\Relbar}
\makeatother

\newcommand{\xeq}{\xequal}

\newcommand{\mto}{\mapsto}

\newcommand{\Set}[2]{\ensuremath{\left\{{#1}\,\middle|\,{#2}\right\}}}

\renewcommand{\k}{\ensuremath{\ol{\mathrm{P}}}}

\newcommand{\h}{\hline}

\newcommand{\ts}[3]{\left[\,{\di #1}\,\right]^{#2}_{#3}}

\renewcommand{\k}[1]{\ensuremath{\left({#1}\right)}}

\newcommand{\ds}{\dots}

\newcommand{\bca}{\begin{cases}}
\newcommand{\eca}{\end{cases}}

\newcommand{\mug}{\ensuremath{\infty}}
\newcommand{\der}[1]{\ensuremath{\lef({#1}\ri)'}}

\newcommand{\cosx}{\ensuremath{\cos x}}

\newcommand{\logx}{\ensuremath{\log x}}

\newcommand{\ff}[2]{\ensuremath{\di\fr{#1}{#2}}}

\renewcommand{\ss}[3]{\ensuremath{\di\int_{#1}^{#2}{#3}\,dx}}

\newcommand{\bpic}{\begin{picture}}\newcommand{\epic}{\end{picture}}

\newcommand{\beda}{\begin{edaenumerate}}
\newcommand{\eeda}{\end{edaenumerate}}

\newcommand{\cd}{\cdots}

\newcommand{\asx}{\ensuremath{\sin^{-1}x}}

\newcommand{\acx}{\ensuremath{\cos^{-1}x}}
\newcommand{\atx}{\ensuremath{\tan^{-1}x}}

\newcommand{\q}{\quad}

\newcommand{\bq}{\begin{quote}}\newcommand{\eq}{\end{quote}}

\newcommand{\rt}{\sqrt}
\newcommand{\be}{\begin{enumerate}}\newcommand{\ee}{\end{enumerate}}
\newcommand{\bce}{\begin{center}}\newcommand{\ece}{\end{center}}
\newcommand{\bde}{\begin{description}}\newcommand{\ede}{\end{description}}
\newcommand{\bri}{\begin{flushright}}\newcommand{\eri}{\end{flushright}}
\newcommand{\bb}{\begin{block}}\newcommand{\eb}{\end{block}}
\newcommand{\bt}{\begin{thm}}\newcommand{\et}{\end{thm}}
\newcommand{\bpf}{\begin{proof}}\newcommand{\epf}{\end{proof}}
\newcommand{\bex}{\begin{ex}}\newcommand{\eex}{\end{ex}}
\newcommand{\bexr}{\begin{exr}}\newcommand{\eexr}{\end{exr}}
\newcommand{\bft}{\begin{fact}}\newcommand{\eft}{\end{fact}}
\newcommand{\brk}{\begin{rmk}}\newcommand{\erk}{\end{rmk}}
\newcommand{\ba}{\begin{align*}}\newcommand{\ea}{\end{align*}}
\newcommand{\bexe}{\begin{exe}}\newcommand{\eexe}{\end{exe}}

\newcommand{\bit}{\begin{itemize}}\newcommand{\eit}{\end{itemize}}

\newcommand{\bcm}{}

\newcommand{\hf}{\hfill}
\newcommand{\fr}{\frac}
\newcommand{\cc}{\ensuremath{\mathbf{C}}}

\newcommand{\rr}{\ensuremath{\mathbf{R}}}

\newcommand{\bd}{\begin{defn}}\newcommand{\ed}{\end{defn}}
\newcommand{\bp}{\begin{prop}}\newcommand{\ep}{\end{prop}}
\newcommand{\p}{\ensuremath{\pi}}
\newcommand{\eh}{\emph}\newcommand{\al}{\alpha}

\newcommand{\te}{\text}\newcommand{\ph}{\phantom}
\newcommand{\wt}{\widetilde}
\newcommand{\lef}{\left}\newcommand{\ri}{\right}

\newcommand{\di}{\displaystyle}

\newcommand{\f}{\frac}

\newcommand{\z}{\ensuremath{\bm{z}}}
\newcommand{\np}{\newpage}

\usepackage[dvipdfmx]{graphicx}
\usepackage[dvipsnames]{xcolor}
\usepackage{float,afterpage}
\usepackage{asymptote,layout}
\usepackage{wrapfig,epic}

\usepackage{geometry}
\usepackage{exscale,latexsym,bm}
\usepackage{amssymb,enumerate,amsmath,amsthm,amsfonts}
\usepackage{verbatim,fancybox,wasysym,fancyhdr,type1cm}
\usepackage[frame,all,poly,curve,knot,arrow]{xy}
\usepackage{colortbl}
\usepackage{boxedminipage}
\usepackage{multirow}

\graphicspath{{./figs/}}
\theoremstyle{definition}
\newtheorem{thm}{Theorem}[section]

\newtheorem{lem}[thm]{Lemma}
\newtheorem{prop}[thm]{Proposition}\newtheorem{cor}[thm]{Corollary}

\newtheorem{exr}[thm]{Exercise}
\newtheorem{ob}[thm]{Observation}

\newtheorem{ex}[thm]{Example}

\newtheorem{defn}[thm]{Definition}\newtheorem{rmk}[thm]{Remark}
\newtheorem{fact}[thm]{Fact}
\newtheorem{block}[thm]{}
\newtheorem*{exe}{Exercise}

\renewcommand{\z}{\zeta}

\begin{document}
\renewcommand{\h}{\hline}
\renewcommand{\arraystretch}{1.5}
\newcommand{\fib}{\text{fib\,}}
\newcommand{\fdp}{\text{FDP\,}}

\title[integral representations for dilogarithm and trilogarithm]{Integral representations for local dilogarithm and trilogarithm functions}
\author{Masato Kobayashi}
\date{\today}                                       
\address{Masato Kobayashi\\
Department of Engineering\\
Kanagawa University, 3-27-1 Rokkaku-bashi, Yokohama 221-8686, Japan.}
\email{masato210@gmail.com}


\maketitle
\begin{abstract}
We show new integral representations for dilogarithm and trilogarithm functions on the unit interval.
As a consequence, we also prove (1) new integral representations for Ap\'{e}ry, Catalan constants and  Legendre $\chi$ functions of order 2, 3, (2) a lower bound for the dilogarithm function on the unit interval,  (3) new Euler sums.
\end{abstract}

\newcommand{\dil}{\text{Li}_{2}}
\renewcommand{\tri}{\text{Li}_{3}}

\renewcommand{\asx}{\sin^{-1}x}
\renewcommand{\acx}{\cos^{-1}x}
\renewcommand{\ph}{\phi}

\section{Introduction}
\label{sec1}
\subsection*{Polylogarithm function}
\label{sub2}
The \eh{polylogarithm function}
\[
\Li_{s}(z)=
\dsum_{n=1}^{\mug}\ff{z^{n}}{n^{s}}
=
z+\ff{z^{2}}{2^{s}}+\ff{z^{3}}{3^{s}}+\cd, \q 
s, z\in\cc, |z|<1
\]
plays a significant role in many areas of number theory; 
its origin, the dilogarithm $\Li_{2}(z)$, dates back to Abel, Euler, Kummer, Landen and Spence etc.
See Kirillov \cite{kirillov}, Lewin \cite{lewin}, Zagier \cite{zagier} for more details.
The main theme of this article is to better 
understand the relation between the dilogarithm,  trilogarithm $\Li_{3}(z)$ functions 
and zeta values $\z(2)$, $\z(3)$ (Ap\'{e}ry constant), 
$\z(4)$ in terms of new integral representations.


%

\subsection*{Main results}
\label{sub2}
First, we wish to briefly explain work of Boo Rim Choe (1987) \cite{boo}, Ewell (1990) \cite{ewell} and Williams-Yue (1993) \cite[p.1582-1583]{yue} which motivated us. Their common idea is that, from Maclaurin series involving $\asx$, they each derived certain infinite sums related to $\z(2)$ and $\z(3)$ with termwise Wallis integral. Figure \ref{f1} gives summary of this.
\begin{table}
\caption{Summary of Boo Rim Choe, Ewell and Williams-Yue's work}
\label{f1}
\begin{center}
\resizebox{.8\width}{!}{
{\renewcommand{\arraystretch}{4}
\begin{tabular}{@{} c|ccc @{}}
\h
Boo Rim Choe&
$\asx=
\dsum_{n=0}^{\mug}\ff{(2n-1)!!}{(2n)!!}\ff{x^{2n+1}}{2n+1}$
&$\to$&
$\dsum_{n=0}^{\mug}\ff{1}{(2n+1)^{2}}=\ff{3}{4}\z(2)=\ff{\p^{2}}{8}$
  \\\hline
Ewell&
$\ff{\asx}{x}=
\dsum_{n=0}^{\mug}\ff{(2n-1)!!}{(2n)!!}\ff{x^{2n}}{2n+1}$
&$\to$&
$\dsum_{n=0}^{\mug}\ff{1}{(2n+1)^{3}}=
\ff{7}{8}\z(3)$
  \\\hline
Williams-Yue&
$\ff{(\asx)^{2}}{x}=
\ff{1}{2}\dsum_{n=1}^{\mug}\ff{(2n)!!}{(2n-1)!!}\ff{x^{2n-1}}{n^{2}}$
&$\to$&
$\ff{\p}{8}
\dsum_{n=1}^{\mug}\ff{1}{n^{3}}=
\ff{\p}{8}\z(3)$
  \\\hline
\end{tabular}
}
}
\end{center}
\end{table}

In this article, we reformulate their ideas introducing  \eh{Wallis operator} and naturally extend their results.
\begin{itemize}
	\item We find new integral representations for 
	$\Li_{2}(t)$, $\Li_{3}(t)$, Legendre $\chi$ functions of order 2, 3 and even for Ap\'{e}ry, Catalan constants  (Theorems \ref{t1}, \ref{t2}, Corollaries \ref{c1}, \ref{c2}). 
	\item We give a lower bound for $\Li_{2}(t)$ on the unit interval (Theorem \ref{t3}).
	\item Making use of $(\sin^{-1}x)^{3}$ and $(\sin^{-1}x)^{4}$, we prove new Euler sums (Theorem \ref{t4}).
\end{itemize}

\subsection*{Notation}
\label{sub3}




Throughout $n$ denotes a nonnegative integer.
Let 
\begin{align*}
	(2n)!!&=2n(2n-2)\cd 4\cdot 2,
	\\(2n-1)!!&=(2n-1)(2n-3)\cd 3\cdot 1.
\end{align*}
In particular, we understand that $(-1)!!=0!!=1$.
Moreover, let \[
w_{n}=\ff{(n-1)!!}{n!!}.
\]
Notice the relation $w_{2n}w_{2n+1}=\tff{1}{2n+1}$ as we will see in the sequel.
\begin{rmk}\hf
\begin{enumerate}
	\item The sequence $\{w_{n}\}$ appears in Wallis integral as 
\[
\ss{0}{\pi/2}{\sin^n x}
=
\begin{cases}
	\f{\p}{2} w_{n}&	\text{$n$ even,}\\
	w_{n}&	\text{$n$ odd.}\\
\end{cases}
\]
\item It also appears in the literature in the disguise of central binomial coefficients as 
\[
w_{2n}=
\ff{(2n-1)!!}{(2n)!!}=\ff{1}{2^{2n}}\binom{2n}{n}.
\]
See Ap\'{e}ry \cite{apery}, van der Poorten \cite{van},  for example.
\end{enumerate}
\end{rmk}

Unless otherwise specified, $t, u, x, y$ are real numbers. 
By $\sin^{-1} x$ and $\cos^{-1} x$, we mean the real inverse sine and cosine functions ($\arcsin x, \arccos x$), that is, 
\[
\begin{array}{ccl}
	y=\asx&\iff   &   x=\sin y, \q-\f{\,\pi\,}{2}\le y\le \f{\,\pi\,}{2},\\
	y=\acx&\iff   &    x=\cos y, \q0\le y\le \p.
\end{array}
\]
\begin{fact}[{Gradshteyn-Ryzhik \cite[p.60, 61]{gr}}]
\begin{eqnarray}
\sin^{-1}t=
\dsum_{n=0}^{\mug}w_{2n}\ff{t^{2n+1}}{2n+1}, \q |t|\le 1.
\label{1}
\end{eqnarray}
\begin{eqnarray}
(\sin^{-1}t)^{2}=
\f12
\dsum_{n=1}^{\mug}
\ff{1}{w_{2n}}\ff{t^{2n}}{n^{2}}, \q |t|\le 1.
\label{2}
\end{eqnarray}
\end{fact}
Further, 
$\sinh^{-1}x=\log(x+\rt{x^{2}+1})$ $(x\inrr)$ 
denotes the inverse hyperbolic sine function 
(some authors write $\te{arsinh\,} x,\te{arcsinh\,} x$ or $\te{argsinh\,} x$ for this one).

\section{Dilogarithm function}
\label{sec2}
\subsection{Definition}

\begin{defn}
For $0\le t\le 1$, the \eh{dilogarithm function} is 
\[
\dil(t)=
\dsum_{n=1}^{\mug}\ff{t^{n}}{n^{2}}.
\]
\end{defn}
In particular, 
$\dil(1)=\z(2)=\tff{\p^{2}}{6}.$ 


It is possible to describe its \eh{even part} by $\dil$ itself since 
\[
\dsum_{n=1}^{\mug}\ff{t^{2n}}{(2n)^{2}}
=\ff{1}{4}\dsum_{n=1}^{\mug}\ff{(t^2)^{n}}{n^{2}}
=\ff{1}{4}\dil(t^{2}).
\]
Its \eh{odd part} is called 
the Legendre $\chi$ function of order 2:
\[
\chi_2(t)=\dsum_{n=1}^{\mug}\ff{t^{2n-1}}{(2n-1)^{2}}.
\]
Here is a fundamental relation of these two parts.
\begin{ob}
\label{ob1}
\[
\dil(t)=\chi_2(t)+\ff{1}{4}\dil(t^{2}).
\]
\end{ob}


\begin{defn}
Define
\[
\te{Ti}_{2}(t)= 
\dsum_{n=1}^{\mug}\ff{(-1)^{n-1}}{(2n-1)^{2}}t^{2n-1}
\]
as a signed analog of $\chi_2(t)$.
\end{defn}
This is also called the 
\eh{inverse tangent integral} of order 2 because of the integral representations 
\[
\te{Ti}_{2}(t)=
\di\int_{0}^{t}{\ff{\tan^{-1}x}{x}}\,dx.
\]

\subsection{Wallis operator}

Let $\rr[[t]]$ denote the set of power series in $t$ 
over real coefficients.
Set \[
F(t)=
\Set{f\in \rr[[t]]}{
\te{$f(t)$ is convergent for $|t|\le 1$} 
}.
\]
\begin{defn}
For $f\in F(t)$, define $W:F(t)\to F(t)$ by 
\[
Wf(t)=
\wi{f(tu)}.
\]
Call $W$ the \eh{Wallis operator}.
\end{defn}

\begin{rmk}
\cite[p.17]{gr} Power series may be integrated and differentiated termwise inside the circle of convergence 
without changing the radius of convergence. 
In the sequel, we will frequently use this without mentioning explicitly.
\end{rmk}	

It is now helpful to understand $W$ coefficientwise.
\begin{lem}
\label{l1}Let $f(t)= \sum_{n=0}^{\mug}a_{n}t^{n}\in F(t)$. 
Then 
\[
Wf(t)=
\dsum_{n=0}^{\mug}a_{2n}\k{\ff{\p}{2}w_{2n}}
t^{2n}+
\dsum_{n=0}^{\mug}a_{2n+1}w_{2n+1}t^{2n+1}.
\]
\end{lem}

\begin{proof}
\begin{align*}
	Wf(t)&=\di\int_{0}^{1}{f(tu)}\,\ff{du}{\rt{1-u^{2}}}
	\\&=\di\int_{0}^{1}{
\k{\dsum_{n=0}^{\mug}a_{2n}t^{2n}u^{2n}
+
\dsum_{n=0}^{\mug}a_{2n+1}t^{2n+1}u^{2n+1}}
}\,\ff{du}{\rt{1-u^{2}}}
	\\&=
\dsum_{n=0}^{\mug}a_{2n}t^{2n}
\di\int_{0}^{1}{u^{2n}}\,\ff{du}{\rt{1-u^{2}}}	
+
\dsum_{n=0}^{\mug}a_{2n+1}t^{2n+1}
\di\int_{0}^{1}{u^{2n+1}}\,\ff{du}{\rt{1-u^{2}}}	
	\\&=
	\dsum_{n=0}^{\mug}a_{2n}\k{\ff{\p}{2}w_{2n}}
t^{2n}+
\dsum_{n=0}^{\mug}a_{2n+1}w_{2n+1}t^{2n+1}.
\end{align*}
\end{proof}
Observe that 
$W$ is linear in the sense that 
$W(f+g)=W(f)+W(g)$ and 
$W(cf)=cW(f)$ for $f, g\in F(t), c\inrr$.


\kill{\np}
\subsection{Main theorem 1}

\begin{lem}
\label{l2}
All of the following are convergent power series for $|t|\le1$.
\begin{eqnarray}
{\sin^{-1}t}=
\dsum_{n=0}^{\mug}w_{2n}\ff{t^{2n+1}}{2n+1}.
\label{3}
\end{eqnarray}
\begin{eqnarray}
\f{\,1\,}{2}(\sin^{-1}t)^{2}=
\dsum_{n=1}^{\mug}\ff{1}{w_{2n}}\ff{t^{2n}}{(2n)^{2}}.
\label{4}
\end{eqnarray}
\begin{eqnarray}
\sin^{-1}t+\ff{1}{\p}(\sin^{-1}t)^{2}
=
\dsum_{n=0}^{\mug}w_{2n}
\ff{t^{2n+1}}{2n+1}+
\dsum_{n=1}^{\mug}\ff{2}{\pi w_{2n}}\ff{t^{2n}}{(2n)^{2}}.
\label{5}
\end{eqnarray}
\begin{eqnarray}
	\sinh^{-1}t=
	\dsum_{n=0}^{\mug}(-1)^{n}
	w_{2n}\ff{t^{2n+1}}{2n+1}.
\label{6}
\end{eqnarray}
\begin{eqnarray}
\f{\,1\,}{2}
(\sinh^{-1}t)^{2}=
\dsum_{n=1}^{\mug}\ff{(-1)^{n-1}}{w_{2n}}\ff{t^{2n}}{(2n)^{2}}.
\label{7}
\end{eqnarray}
\end{lem}

\begin{proof}
We already saw (\ref{3}) and (\ref{4}) in  Introduction.
(\ref{5}) is $(\ref{3})+\tf{2}{\p}(\ref{4})$.
(\ref{6}) and (\ref{7}) follow from 
(\ref{3}), (\ref{4}) 
and $\sinh^{-1}z=\tf{1}{i}\sin^{-1}(iz)$ (for all $z\in\cc$) 
\cite[p.56]{gr}.
\end{proof}


\begin{thm}\label{t1}
For $0\le t\le 1$, all of the following hold.
\begin{eqnarray}
\chi_2(t)=
\di\int_{0}^{1}{\ff{\sin^{-1}(tu)}{\rt{1-u^{2}}}}\,du.
\label{8}
\end{eqnarray}
\begin{eqnarray}
\ff{1}{4}\Li_{2}(t^{2})=
\ff{1}{\p} \di\int_{0}^{1}{\ff{(\sin^{-1}(tu))^{2}}{\rt{1-u^{2}}}}\,du.
\label{9}
\end{eqnarray}
\begin{eqnarray}
\Li_{2}(t)=
\di\int_{0}^{1}{\ff{\sin^{-1}(tu)+\tf{1}{\p}(\sin^{-1}(tu))^{2}}{\rt{1-u^{2}}}}\,du.
\label{10}
\end{eqnarray}
\begin{eqnarray}
\te{Ti}_{2}(t)=
\di\int_{0}^{1}{\ff{\sinh^{-1}(tu)}{\rt{1-u^{2}}}}\,du.
\label{11}
\end{eqnarray}
\begin{eqnarray}
\ff{\pi}{2} 
\k{\ff{1}{4}\Li_{2}(t^{2})-\ff{\,1\,}{8}\Li_{2}(t^{4})}=
\di\int_{0}^{1}{
\ff{\tf12(\sinh^{-1}tu)^{2}}{\rt{1-u^{2}}}
}\,du
.
\label{12}
\end{eqnarray}
\end{thm}

\begin{proof}
Note that these are equivalent to the following statements:
\begin{eqnarray}
W\k{{\sin^{-1}t}}
=\chi_2(t).
\label{13}
\end{eqnarray}
\begin{eqnarray}
W
\k{
{\f{\,1\,}{2}
(\sin^{-1}t)^{2}}}
=
\ff{\p}{2}\cdot \ff{1}{4}\te{Li}_{2}(t^{2}).
\label{14}
\end{eqnarray}
\begin{eqnarray}
W
\k{\sin^{-1}t+\ff{1}{\p}(\sin^{-1}t)^{2}}
=
\te{Li}_{2}(t).
\label{15}
\end{eqnarray}
\begin{eqnarray}
	W(\sinh^{-1}t)=\te{Ti}_{2}(t).
\label{16}
\end{eqnarray}
\begin{eqnarray}
W\k{
{\f{\,1\,}{2}
(\sinh^{-1}t)^{2}}}
=
\ff{\pi}{2} 
\k{\ff{1}{4}\Li_{2}(t^{2})-\ff{\,1\,}{8}\Li_{2}(t^{4})}.
\label{17}
\end{eqnarray}
With Lemmas \ref{l1} and \ref{l2}, we can verify (\ref{13})-(\ref{16}) by checking  coefficients of those series. For example, 
\begin{align*}
	W(\sin^{-1}t)&=W\k{
\dsum_{n=0}^{\mug}w_{2n}\ff{t^{2n+1}}{2n+1}
}
	=\dsum_{n=0}^{\mug}w_{2n}w_{2n+1}\ff{t^{2n+1}}{2n+1}
\\&=\dsum_{n=0}^{\mug}\ff{t^{2n+1}}{(2n+1)^{2}}=\chi_2(t).
\end{align*}
It remains to show (\ref{17}).
\begin{align*}
	W\k{
{\f{\,1\,}{2}
(\sinh^{-1}t)^{2}}}
&=\dsum_{n=1}^{\mug}\ff{(-1)^{n-1}}{w_{2n}}
\k{w_{2n}\ff{\pi }{2}}
\ff{t^{2n}}{(2n)^{2}}
	\\&=\ff{\pi }{2}
	\dsum_{n=1}^{\mug}\ff{(-1)^{n-1}}{(2n)^{2}}
	t^{2n}
	\\&=\ff{\pi }{2}
	\k{
	\dsum_{n=1}^{\mug}\ff{1}{(2n)^{2}}
	t^{2n}
	-
	2
	\dsum_{n=1}^{\mug}\ff{1}{(4n)^{2}}
	t^{4n}
	}
	\\&=
	\ff{\pi}{2} 
\k{\ff{1}{4}\Li_{2}(t^{2})-\ff{\,1\,}{8}\Li_{2}(t^{4})}.
\end{align*}
\end{proof}

\kill{\np}

\begin{cor}
\label{c1}
\begin{eqnarray}
\di\int_{0}^{1}{}\,\ff{\sin^{-1}u }{\rt{1-u^{2}}}du
=\f{3}{4}\z(2)=\f{\p^{2}}{8}.
\label{}
\end{eqnarray}
\begin{eqnarray}
\ff{2}{\pi}
\di\int_{0}^{1}{}\,\ff{\tf12(\sin^{-1}u)^{2}}{\rt{1-u^{2}}}du
=\ff{1}{4}\z(2)=
\ff{\p^{2}}{24}.
\label{}
\end{eqnarray}	
\begin{eqnarray}
	\di\int_{0}^{1}{
\k{\sin^{-1}u+\ff{1}{\pi}(\sin^{-1}u)^{2}}
}\,\ff{du}{\rt{1-u^{2}}}
=\z(2)=\ff{\p^{2}}{6}.
\label{}
\end{eqnarray}
\begin{eqnarray}
\di\int_{0}^{1}{\ff{\sinh^{-1}u}{\rt{1-u^{2}}}}\,du=G.
\label{}
\end{eqnarray}
\begin{eqnarray}
	\di\int_{0}^{1}{\ff{\tf{1}{2}(\sinh^{-1}u)^{2}}{\rt{1-u^{2}}}}\,du=\ff{\p}{16}\z(2)=\ff{\p^{3}}{96}.
\label{}
\end{eqnarray}
\end{cor}

\begin{proof}
These are 
$\chi_{2}(1), \tf{1}{4} \dil(1^{2}), \dil(1), \te{Ti}_{2}(1)$ and 
$\tff{\pi}{2} 
\k{\tff{1}{4}\Li_{2}(1^{2})-\tff{\,1\,}{8}\Li_{2}(1^{4})}$.
\end{proof}

\kill{\np}

\section{Trilogarithm function}
\label{sec3}
\subsection{Definition}

\begin{defn}The \eh{trilogarithm function} for $0\le t\le 1$ is 
\[
\te{Li}_{3}(t)=
\dsum_{n=1}^{\mug}\ff{t^{n}}{n^{3}}.
\]
\end{defn}
Its odd part is the Legendre $\chi$ function of order 3:
\[
\chi_3(t)=
\dsum_{n=1}^{\mug}\ff{t^{2n-1}}{(2n-1)^{3}}.
\]
In particular, $\te{Li}_{3}(1)=\z(3)$ and 
$\chi_{3}(1)=\tf{7}{8}\z(3)$.

\begin{ob}
\label{ob2}
\[
\tri(t)=\chi_3(t)+\ff{1}{8}\tri(t^{2}).
\]
\end{ob}

Further, a signed analog of $\chi_3(t)$ is 
\[
\te{Ti}_{3}(t)=
\dsum_{n=1}^{\mug}\ff{(-1)^{n-1}}{(2n-1)^{3}}t^{2n-1}.
\]

\subsection{Main theorem 2}

\begin{lem}
\label{l3}
\begin{eqnarray}
\di\int_{0}^{t}{\ff{\sin^{-1}y}{y}}\,dy=
\dsum_{n=0}^{\mug}w_{2n}\ff{t^{2n+1}}{(2n+1)^{2}}.
\label{23}
\end{eqnarray}
\begin{eqnarray}
\di\int_{0}^{t}{\ff{\tf12(\sin^{-1}y)^{2}}{y}}\,dy=
\dsum_{n=1}^{\mug}\ff{1}{w_{2n}}\ff{t^{2n}}{(2n)^{3}}.
\label{24}
\end{eqnarray}
\begin{eqnarray}
\di\int_{0}^{t}{\ff{\sin^{-1}y+\f{1}{\p}(\sin^{-1}y)^{2}}{y}}\,dy
=
\dsum_{n=0}^{\mug}w_{2n}\ff{t^{2n+1}}{(2n+1)^{2}}
+
\dsum_{n=1}^{\mug}\ff{2}{\pi w_{2n}}\ff{t^{2n}}{(2n)^{3}}.
\label{25}
\end{eqnarray}
\begin{eqnarray}
\di\int_{0}^{t}{\ff{\sinh^{-1}y}{y}}\,dy
=
\dsum_{n=0}^{\mug}(-1)^{n}w_{2n}\ff{t^{2n+1}}{(2n+1)^{2}}.
\label{26}
\end{eqnarray}
\begin{eqnarray}
\di\int_{0}^{t}{\ff{\tf12(\sinh^{-1}y)^{2}}{y}}\,dy=
\dsum_{n=1}^{\mug}\ff{(-1)^{n-1}}{w_{2n}}\ff{t^{2n}}{(2n)^{3}}.
\label{27}
\end{eqnarray}
\end{lem}

\begin{proof}
We can derive all of these 
by integrating (\ref{3})-(\ref{7}) termwise.
\end{proof}
As a consequence, we obtain the 
equalities below (cf. (\ref{13})-(\ref{17})).
\begin{eqnarray}
W\k{
\di\int_{0}^{t}{\ff{\sin^{-1}y}{y}}\,dy}=\chi_3(t).
\label{28}
\end{eqnarray}
\begin{eqnarray}
W\k{
\di\int_{0}^{t}{\ff{\tf{1}{2}(\sin^{-1}y)^{2}}{y}}\,dy}=\ff{\p}{2}\cdot \ff{1}{8}\tri(t^{2}).
\label{29}
\end{eqnarray}
\begin{eqnarray}
W\k{
\di\int_{0}^{t}{\ff{\sin^{-1}y+\f{1}{\p}(\sin^{-1}y)^{2}}{y}}\,dy}=\tri(t).
\label{30}
\end{eqnarray}
\begin{eqnarray}
W\k{
\di\int_{0}^{t}{\ff{\sinh^{-1}y}{y}}\,dy}=\te{Ti}_{3}(t).
\label{31}
\end{eqnarray}
\begin{eqnarray}
W\k{
\di\int_{0}^{t}{\ff{\tf{1}{2}(\sinh^{-1}y)^{2}}{y}}\,dy}
=\dsum_{n=1}^{\mug}\ff{(-1)^{n-1}}{w_{2n}}
\k{w_{2n}\ff{\pi }{2}}
\ff{t^{2n}}{(2n)^{3}}
\label{32}
\end{eqnarray}
\begin{align*}
	\phantom{W\k{
\di\int_{0}^{t}{\ff{\tf{1}{2}(\sinh^{-1}y)^{2}}{y}}\,dy}}
	&=\ff{\pi }{2}
	\dsum_{n=1}^{\mug}\ff{(-1)^{n-1}}{(2n)^{3}}
	t^{2n}
	\\&=\ff{\pi }{2}
	\k{
	\dsum_{n=1}^{\mug}\ff{1}{(2n)^{3}}
	t^{2n}
	-
	2
	\dsum_{n=1}^{\mug}\ff{1}{(4n)^{3}}
	t^{4n}
	}
	\\&=\ff{\pi}{2} 
\k{\ff{1}{8}\Li_{3}(t^{2})-\ff{\,1\,}{32}\Li_{3}(t^{4})}.
\end{align*}

In this way, the five functions above come to possess \eh{double} integral representations.
For example, 
\[
\chi_3(t)=
\di\int_{0}^{1}{
\k{
\di\int_{0}^{tu}{\ff{\sin^{-1}y}{y}}\,dy
}
}\,\ff{du}{\rt{1-u^{2}}}.
\]
We can indeed simplify such integrals to \eh{single} ones by exchanging order of integrals.
\begin{thm}
\label{t2}
\begin{eqnarray}
\chi_3(t)=
\di\int_{0}^{1}{\ff{\sin^{-1}(tx)\cos^{-1}x}{x}}\,dx.
\label{33}
\end{eqnarray}
\begin{eqnarray}
\ff{1}{8}\te{Li}_{3}(t^{2})=
\ff{2}{\pi}
\di\int_{0}^{1}{
\ff{\tf{1}{2}(\sin^{-1}(tx))^{2}\acx}{x}
}\,dx.
\label{34}
\end{eqnarray}
\begin{eqnarray}
\te{Li}_{3}(t)=
\di\int_{0}^{1}{
\ff{
\k{\sin^{-1}(tx)+\tf{1}{\p}(\sin^{-1}(tx))^{2}}
\cos^{-1}x}
{x}
}\,dx.	
\label{35}
\end{eqnarray}
\begin{eqnarray}
\te{Ti}_{3}(t)=\di\int_{0}^{1}\ff{\sinh^{-1}(tx)\acx}{x}dx.
\label{36}
\end{eqnarray}
\begin{eqnarray}
\ff{\p}{2}
\k{\ff{1}{8}\tri(t^{2})-\ff{1}{32}\tri(t^{4})}
=\di\int_{0}^{1}\ff{\tf12(\sinh^{-1}(tx))^{2}\acx}{x}dx.
\label{37}
\end{eqnarray}
\end{thm}

\begin{proof}
We give a proof altogether. 
For $t=0$, all the equalities hold as $0=0$. 
Suppose $0<t\le 1$. Let 
\[
f(y)\in\left\{
\sin^{-1}y, 
\ff{1}{\p}(\sin^{-1}y)^{2}, 
\sin^{-1}y+\ff{1}{\p}(\sin^{-1}y)^{2},
\sinh^{-1}y, 
\ff{1}{2}(\sinh^{-1}y)^{2}\right\}.
\]
Then 
\begin{align*}
	W\k{
\di\int_{0}^{t}{\ff{f(y)}{y}}\,dy}
&=
\di\int_{0}^{1}{
\di\int_{0}^{tu}{\ff{f(y)}{y}}\,dy
}\,\ff{du}{\rt{1-u^{2}}}
	\\&=
\di\int_{0}^{t}{
\di\int_{y/t}^{1}{\ff{f(y)}{y}}
\ff{1}{\rt{1-u^{2}}}
}\,dudy
	\\&=\di\int_{0}^{t}
	{\ff{f(y)}{y}}\cos^{-1}{\ff{y}{t}} \,dy
	\\&=\di\int_{0}^{1}
	{\ff{f(tx)}{x}}\cos^{-1}{x} \,dx.
\end{align*}

\end{proof}

\kill{\np}

\begin{cor}
\label{c2}
\begin{eqnarray}
\di\int_{0}^{1}\ff{\sin^{-1}x\cos^{-1}x}{x}dx
=\ff{7}{8}\z(3).
\label{38}
\end{eqnarray}
\begin{eqnarray}
\ff{2}{\p}\di\int_{0}^{1}\ff{\tf{1}{2}(\sin^{-1}x)^{2}\cos^{-1}x}{x}dx
=\ff{1}{8}\z(3)
.
\label{39}
\end{eqnarray}
\begin{eqnarray}
\di\int_{0}^{1}\ff{
\k{\sin^{-1}x+\tf{1}{\p}(\sin^{-1}x)^{2}}
\cos^{-1}x}{x}
dx=\z(3)
.
\label{40}
\end{eqnarray}
\begin{eqnarray}
\di\int_{0}^{1}\ff{\sinh^{-1}x\acx}{x}dx
=\ff{\p^{3}}{32}
.
\label{41}
\end{eqnarray}
\begin{eqnarray}
\di\int_{0}^{1}\ff{\tf12(\sinh^{-1}x)^{2}\cos^{-1}x}{x}dx
=\ff{3\p}{64}\z(3)
.
\label{42}
\end{eqnarray}

\end{cor}

\begin{proof}
These are $\chi_{3}(1), 
\tf{1}{8}\tri(1^{2}), \tri(1), \te{Ti}_{3}(1)$ 
 and 
 $\tff{\p}{2}\k{\tff{1}{8}\tri(1^{2})-\tff{1}{32}\tri(1^{4})}$.
\end{proof}


\kill{\np}

\kill{\np}

\kill{\np}




\kill{\np}

\section{Applications}
\label{sec4}

\subsection{Inequalities}

It is easy to see from 
the definitions 
$\dil(t)=\sum_{n=1}^{\mug}\tf{t^{n}}{n^{2}}$ 
and 
$\chi_2(t)=\sum_{n=1}^{\mug}\tfrac{t^{2n-1}}{(2n-1)^{2}}$ ($0\le t\le 1$) 
that 
\[
0\le \dil(t)\le \f{\p^{2}}{6}
\qte{and}\q
0\le \chi_2(t)\le \f{\p^{2}}{8}.
\]
In fact, we can improve these inequalities a little more. 
For upper bounds, it is immediate that 
\[
\dil(t)=
\dsum_{n=1}^{\mug}\ff{t^{n}}{n^{2}}
\le
\dsum_{n=1}^{\mug}\ff{t}{n^{2}}
=\ff{\p^{2}}{6}t,
\]
\[
\chi_{2}(t)=
\dsum_{n=1}^{\mug}\ff{t^{2n-1}}{(2n-1)^{2}}
\le
\dsum_{n=1}^{\mug}\ff{t}{(2n-1)^{2}}
=\ff{\p^{2}}{8}t.
\]
We next prove nontrivial lower bounds for these functions and also $\te{Ti}_{2}(t)$.
\begin{thm}\label{t3}
For $0\le t\le 1$,
\begin{eqnarray}
\dil(t)\ge
\ff{4}{3\p}\ff{(\sin^{-1}t)^{3}}{t}
.
\label{43}
\end{eqnarray}
\begin{eqnarray}
\chi_2(t)\ge
\ff{(\sin^{-1}t)^{2}}{2t}
.
\label{44}
\end{eqnarray}
\begin{eqnarray}
\te{Ti}_{2}(t)\ge 
\ff{(\sinh^{-1}t)^{2}}{2t}.
\label{45}
\end{eqnarray}
\end{thm}

%
%

\kill{\np}
Before the proof, we need a lemma. It provides another integral representation of $\dil(t)$ which seems  interesting itself.
\begin{lem}
\label{l4}
For $0\le t\le 1$, 
\begin{eqnarray}
\dil(t)=
\ff{8\rt{t}}{\p}
\di\int_{0}^{1}{
\ff{\sin^{-1}(\rt{t}x)\acx}{\rt{1-tx^{2}}}
}\,dx.
\label{46}
\end{eqnarray}
\[
\k{\te{cf.}
\q 
\tri(t)=\ff{8}{\p}\di\int_{0}^{1}\ff{(\sin^{-1}(\rt{t}x))^{2}\acx}{x}dx,
\q \te{$t\mto \rt{t}$ in (\ref{34}).}
}\]
\end{lem}

\begin{proof}
If $t=0$, then both sides are $0$. For $0<t\le 1$, 
\begin{align*}
	\te{RHS}&=
	\ff{8}{\p} \di\int_{0}^{\rt{t}}{\ff{\sin^{-1}y}{\rt{1-y^{2}}}}\cos^{-1}\ff{y}{\rt{t}}\,dy
	\\&=
	\ff{8}{\p} \di\int_{0}^{\rt{t}}{\ff{\sin^{-1}y}{\rt{1-y^{2}}}}
\di\int_{y/\rt{t}}^{1}{\ff{1}{\rt{1-u^{2}}}}\,dudy
	\\&=\ff{8}{\p} \di\int_{0}^{1}
	\di\int_{0}^{\rt{tu}}
	{\ff{\sin^{-1}y}{\rt{1-y^{2}}}}dy
\ff{1}{\rt{1-u^{2}}}\,du
	\\&=\ff{8}{\p} 
	W\k{\di\int_{0}^{\rt{t}}
	{\ff{\sin^{-1}y}{\rt{1-y^{2}}}}
	dy}
	\\&=\ff{8}{\p} 
	W\k{\f12\k{\sin^{-1}\rt{t}}
^{2}}
	\\&=\ff{8}{\p} 
	\k{\ff{\p}{2} \f14 \dil\k{(\rt{t})^{2}}
}
=\dil(t).
\end{align*}

\end{proof}

\kill{\np}

\begin{proof}[Proof of Theorem \ref{t3}]
If $t=0$, then all of (\ref{43})-(\ref{45}) hold as $0\ge0$.
Suppose $t>0$. 
Since $\sin^{-1}$ is increasing on $[0, 1]$, $\sin^{-1}(tx)\le \sin^{-1}(\rt{t}x)$ for all $0< t, x\le 1$. 
Then 
\begin{align*}
	\dil(t)
&=
\ff{8\rt{t}}{\p}\di\int_{0}^{1}{
\ff{\sin^{-1}(\rt{t}x)}{\rt{1-tx^{2}}}
}\acx\,dx
	\\&\ge
\ff{8t}{\p}\di\int_{0}^{1}{
\ff{\sin^{-1}(tx)}{\rt{1-t^{2}x^{2}}}
}\acx\,dx
	\\&=
\ff{8{t}}{\p}\di\int_{0}^{1}{
\ff{1}{t}\der{\ff{1}{2}(\sin^{-1}tx)^{2}}
}\acx\,dx	
	\\&=
	\ff{8}{\p}
	\k{
	\underbrace{\ts{
	{\ff{1}{2}(\sin^{-1}tx)^{2}}\acx
	}{1}{0}}_{0}
	-
	\di\int_{0}^{1}{\ff{1}{2}(\sin^{-1}tx)^{2}\ff{-1}{\rt{1-x^{2}}}}\,dx
	}
	\\&\ge \ff{4}{\p}
	\di\int_{0}^{1}{
	\ff{(\sin^{-1}tx)^{2}}{\rt{1-t^{2}x^{2}}}}
	\,dx
	\\&=\ff{4}{\p}\ts{\ff{1}{3t}(\sin^{-1}tx)^{3}}{1}{0}
	\\&=\ff{4}{3\p}\ff{(\sin^{-1}t)^{3}}{t}.
\end{align*}
Next, we prove (\ref{44}). Note that 
\[
\ff{\sin^{-1}(tx)}{\rt{1-t^{2}x^{2}}}\le 
\ff{\sin^{-1}(tx)}{\rt{1-x^{2}}}
\]
for $0< t, x< 1$.
Integrate these from $0$ to $1$ in $x$ so that 
\[
\di\int_{0}^{1}\ff{\sin^{-1}(tx)}{\rt{1-t^{2}x^{2}}}
\,dx
\le 
\di\int_{0}^{1}
\ff{\sin^{-1}(tx)}{\rt{1-x^{2}}}
\,
dx,
\]
\[
\ts{\ff{(\sin^{-1}(tx))^{2}}{2t}}{1}{0}\le 
\chi_2(t),
\]
\[
\ff{(\sin^{-1}t)^{2}}{2t}
\le 
\chi_2(t).
\]
Quite similarly, 
for $0< t, x< 1$, it also holds that 
\[
\ff{\sinh^{-1}(tx)}{\rt{1+t^{2}x^{2}}}\le 
\ff{\sinh^{-1}(tx)}{\rt{1-x^{2}}},
\]
\[
\di\int_{0}^{1}{
\ff{\sinh^{-1}(tx)}{\rt{1+t^{2}x^{2}}}}\,dx
\le 
\di\int_{0}^{1}{
\ff{\sinh^{-1}(tx)}{\rt{1-x^{2}}}
}\,dx
=
\te{Ti}_{2}(t).
\]
The left hand side is 
\[
\ts{\ff{\sinh^{-1}(tx)^{2}}{2t}}{1}{0}=\ff{(\sinh^{-1}t)^{2}}{2t}.
\]
\end{proof}


\kill{\np}

\kill{\np}

\kill{\np}

\subsection{Euler sums}

\begin{defn}
A \eh{harmonic number} is $H_{n}=\sum_{k=1}^{n} \tf{1}{k}$.
More generally, for $m, n\ge 1$, 
an $(m, n)$-\eh{harmonic number} is
\[
H_{n}^{(m)}=
\dsum_{k=1}^{n}\ff{1}{k^{m}}.
\]
\end{defn}
In particular, $H_{n}^{(1)}=H_{n}$. 
Any series involving such numbers is called an \eh{Euler sum}.

V\u{a}lean \cite[p.292-293]{valean} presents truly remarkable Euler sums such as 
\begin{align*}
	\dsum_{n=1}^{\mug}\ff{H_{n}^{2}}{n^{2}}&=\ff{17}{4}\z(4),
	\\\dsum_{n=1}^{\mug}\ff{H_{n}^{2}}{n^{3}}&=\ff{7}{2}\z(5)-\z(2)\z(3),
	\\\dsum_{n=1}^{\mug}\ff{H_{n}^{2}}{n^{4}}&=
\ff{97}{24}\z(6)-2\z^{2}(3),
	\\\dsum_{n=1}^{\mug}\ff{H_{n}^{2}}{n^{5}}&=6\z(7)-\z(2)\z(5)-\ff{5}{2}\z(3)\z(4),
	\\\dsum_{n=1}^{\mug}\ff{H_{n}H_{n}^{(2)}}{n^{2}}&=\z(2)\z(3)+\z(5).
\end{align*}
There are many ideas to prove such formulas;
 J.M. Borwein and Bradley \cite{bobr} 
gives thirty two proofs for 
\[
\dsum_{n=1}^{\mug}\ff{H_{n-1}}{n^{2}}=
\z(3)=
8\dsum_{n=1}^{\mug}\ff{(-1)^{n}H_{n-1}}{n^{2}}
\]
by integrals, polylogarithm functions, Fourier series and hypergeometric functions etc.
Here, as an application of our main idea, Wallis operators, we prove two new Euler sums. 
Let 
\[
{O}_{n}^{(2)}=H_{2n-1}^{(2)}-\ff{1}{4}H_{n-1}^{(2)}
=
\dsum_{k=1}^{2n-1}\ff{1}{k^{2}}-
\dsum_{k=1}^{n-1}\ff{1}{(2k)^{2}}=
\dsum_{k=0}^{n-1} 
\ff{1}{(2k+1)^{2}}.
\]

\begin{thm}\label{t4}
\begin{eqnarray}
\dsum_{n=0}^{\mug}\ff{{O}_{n}^{(2)}}{(2n+1)^{2}}
=\ff{\p^{4}}{384}=\ff{15}{64}\z(4).
\label{47}
\end{eqnarray}
\begin{eqnarray}
\dsum_{n=1}^{\mug}\ff{{H}_{n-1}^{(2)}}{n^{2}}
=\ff{\p^{4}}{120}=\ff{3}{4}\z(4).
\label{48}
\end{eqnarray}
\end{thm}

For the proof, we make use of less-known 
Maclaurin series for $(\sin^{-1}t)^{3}$ and 
$(\sin^{-1}t)^{4}$; thus we can interpret this result as a natural subsequence of Boo, Ewell and Williams-Yue's work.

\begin{lem}
\label{l5}
\begin{eqnarray}
(\sin^{-1}t)^{3}=
\dsum_{n=0}^{\mug}
\k{
6O_{n}^{(2)}
}
w_{2n}
\ff{t^{2n+1}}{2n+1}.
\label{49}
\end{eqnarray}
\begin{eqnarray}
(\sin^{-1}t)^{4}=
\ff{\,1\,}{2}
\dsum_{n=1}^{\mug}
\k{3H_{n-1}^{(2)}}
\ff{1}{w_{2n}}\ff{t^{2n}}{n^{2}}.
\label{50}
\end{eqnarray}
\[
\k{
\te{cf.}\q 
\sin^{-1}t=
\dsum_{n=0}^{\mug}
w_{2n}
\ff{t^{2n+1}}{2n+1},
\q 
(\sin^{-1}t)^{2}=
\ff{\,1\,}{2}
\dsum_{n=1}^{\mug}
\ff{1}{w_{2n}}\ff{t^{2n}}{n^{2}}
}
\]
\end{lem}

\begin{proof}
First, write 
$(\sin^{-1}t)^{3}=
\sum_{n=0}^{\mug}A_{n}t^{2n+1}$, $A_{n}\inrr$ and 
let $a_{n}=\tff{2n+1}{w_{2n}}A_{n}\, (n\ge0)$. 
It is enough to show that $a_{n}=6O_{n}^{(2)}$. 
Since 
the series $(\sin^{-1}t)^{3}=(t+\tf{t^{3}}{6}+\cd)^{3}$ starts from the $t^{3}$ term, $A_{0}=a_{0}=0$.
For convenience, set 
\[
f_{n}(x)=\ff{\sin^{2n+1}x}{(2n+1)!}.
\]
Then 
\[
f_{n}'(x)=\ff{\sin^{2n}x}{(2n)!}\cosx,
\]
\[
f_{n}''(x)=\ff{1}{(2n)!}
\k{2n\sin^{2n-1}x(1-\sin^{2}x)-\sin^{2n+1}x}
=f_{n-1}(x)-(2n+1)^{2}f_{n}(x).
\]
Now let $x=\sin^{-1}t$ ($-\tf{\,\pi\,}{2}\le x\le \tf{\,\pi\,}{2}$), 
$b_{n}=(2n-1)!!$.
Recall that 
\[
\sin^{-1}t=
\dsum_{n=0}^{\mug}
w_{2n}
\ff{t^{2n+1}}{2n+1}.
\]
In terms of $x$, $b_{n}, f_{n}(x)$, this is 
\[
x=
\dsum_{n=0}^{\mug}
w_{2n}
\ff{\sin^{2n+1}x}{2n+1}
=
\dsum_{n=0}^{\mug}
\ff{(2n-1)!!}{(2n)!!}(2n)!
\ff{\sin^{2n+1}x}{(2n+1)!}
= \dsum_{n=0}^{\mug}b_{n}^{2} f_{n}(x).
\]
Thus, 
\begin{align*}
	x^{3}&=\dsum_{n=0}^{\mug} A_{n}\sin^{2n+1}x
	=\dsum_{n=0}^{\mug} a_{n}
	\k{w_{2n}\ff{\sin^{2n+1}x}{2n+1}}
	=\dsum_{n=0}^{\mug} a_{n}b_{n}^{2}
	f_{n}(x).
\end{align*}
Differentiate both sides twice in $x$:
\begin{align*}
	6x&=
	\dsum_{n=0}^{\mug}a_{n}b_{n}^{2}f_{n}''(x)
	\\&=
\dsum_{n=0}^{\mug}a_{n}b_{n}^{2}(f_{n-1}(x)-(2n+1)^{2}f_{n}(x))
	\\&=
\dsum_{n=0}^{\mug}
(a_{n+1}b_{n+1}^{2}f_{n}(x)-
a_{n}b_{n}^{2}(2n+1)^{2}f_{n}(x)),
\\
6
\dsum_{n=0}^{\mug}b_{n}^{2}f_{n}(x)
&=
\dsum_{n=0}^{\mug}
(a_{n+1}b_{n+1}^{2}f_{n}(x)-
a_{n}b_{n}^{2}(2n+1)^{2}f_{n}(x)).
\end{align*}
Equating coefficients of $f_{n}(x)$ yields
\[
6b_{n}^{2}=a_{n+1}b_{n+1}^{2}-
a_{n}b_{n}^{2}(2n+1)^{2}, \q n\ge0.
\]
Since $b_{n+1}=(2n+1)b_{n}$ and $b_{n}\ne0$, we must have 
\[
a_{n+1}-a_{n}=\ff{6}{(2n+1)^{2}}.
\]
With $a_{0}=0$, we now arrive at 
\[
a_{n}= \dsum_{k=0}^{n-1}\ff{6}{(2k+1)^{2}}=6O_{n}^{(2)},
\]
as required.

The proof for (\ref{50}) proceeds along the same line.
Write 
$(\sin^{-1}t)^{4}=
\tf12\sum_{n=0}^{\mug}C_{n}t^{2n}$, $C_{n}\inrr$ and 
let $c_{n}=C_{n}w_{2n}n^{2}\, (n\ge1)$. 
It is enough to show that $c_{n}=3H_{n-1}^{(2)}$.
Since 
the series $(\sin^{-1}t)^{4}$ starts from the $t^{4}$ term, $C_{1}=c_{1}=0$.
For convenience, set
\[
g_{n}(x)=\ff{\sin^{2n}x}{(2n)!}.
\]
Then 
\[
g_{n}'(x)=\ff{\sin^{2n-1}x}{(2n-1)!}\cosx,
\]
\[
g_{n}''(x)=\ff{1}{(2n-1)!}
\k{(2n-1)\sin^{2n-2}x(1-\sin^{2}x)-\sin^{2n}x}
=g_{n-1}(x)-(2n)^{2}g_{n}(x).
\]
Now let $x=\sin^{-1}t$ ($-\tf{\,\pi\,}{2}\le x\le \tf{\,\pi\,}{2}$), 
$d_{n}=2^{n}(n-1)!$.
Recall that 
\[
(\sin^{-1}t)^{2}=
\ff{1}{2}
\dsum_{n=1}^{\mug}
\ff{1}{w_{2n}}
\ff{t^{2n}}{n^{2}}.
\]
In terms of $x$, $d_{n}, g_{n}(x)$, this is 
\[
x^{2}=
\f12
\dsum_{n=1}^{\mug}
\ff{1}{w_{2n}}
\ff{\sin^{2n}x}{n^{2}}
=
\f12
\dsum_{n=1}^{\mug}
\ff{(2n)!!}{(2n-1)!!}
\ff{(2n)!}{n^{2}}
\ff{\sin^{2n}x}{(2n)!}
= \ff{1}{2}\dsum_{n=1}^{\mug}d_{n}^{2} g_{n}(x).
\]
Thus, 
\begin{align*}
	x^{4}=
	\ff{1}{2}\dsum_{n=1}^{\mug} c_{n}
	\k{\ff{1}{w_{2n}}\ff{\sin^{2n}x}{n^{2}}}=
	\ff{1}{2}\dsum_{n=1}^{\mug} c_{n}
	d_{n}^{2}g_{n}(x).
\end{align*}
Differentiate both sides twice in $x$:
\begin{align*}
	12x^{2}&=
	\f12\dsum_{n=1}^{\mug}c_{n}d_{n}^{2}g_{n}''(x)
	\\&=
	\f12
\dsum_{n=1}^{\mug}c_{n}d_{n}^{2}(g_{n-1}(x)-(2n)^{2}g_{n}(x))
	\\&=
	\f12
\dsum_{n=0}^{\mug}
c_{n+1}d_{n+1}^{2}g_{n}(x)-
\f12
\dsum_{n=1}^{\mug}
c_{n}d_{n}^{2}(2n)^{2}g_{n}(x)
	\\&=
	\f12
\dsum_{n=1}^{\mug}
c_{n+1}d_{n+1}^{2}g_{n}(x)-
\f12
\dsum_{n=1}^{\mug}
c_{n}d_{n}^{2}(2n)^{2}g_{n}(x)
\q\te{($c_{1}=0$)}
\\
&=
\f12
\dsum_{n=1}^{\mug}
(c_{n+1}d_{n+1}^{2}-c_{n}d_{n}^{2}(2n)^{2})g_{n}(x),
\end{align*}
\[
\ff{12}{2}\dsum_{n=1}^{\mug}d_{n}^{2} g_{n}(x)
=
\f12
\dsum_{n=1}^{\mug}
(c_{n+1}d_{n+1}^{2}-c_{n}d_{n}^{2}(2n)^{2})g_{n}(x).
\]
Equating coefficients of $g_{n}(x)$ yields 
\[
\ff{12}{2}d_{n}^{2}=
\f12
(c_{n+1}d_{n+1}^{2}-c_{n}d_{n}^{2}(2n)^{2}), \q n\ge 1.
\]
Since $d_{n+1}=2n d_{n}$ and $d_{n}\ne0$, we must have 
\[
c_{n+1}-c_{n}=\ff{12}{(2n)^{2}}.
\]
With $c_{1}=0$, we conclude that 
\[
c_{n}= \dsum_{k=1}^{n-1}\ff{12}{(2k)^{2}}
=
\dsum_{k=1}^{n-1}\ff{3}{k^{2}}=3H_{n-1}^{(2)}.
\]

\end{proof}

\begin{proof}[Proof of Theorem \ref{t4}]
Note that 
\begin{align*}
	W\k{\ff{1}{6}(\sin^{-1}t)^{3}}&=
\dsum_{n=1}^{\mug}\ff{\wt{O}_{n}^{(2)}w_{2n}}{2n+1}w_{2n+1}t^{2n+1}	
	\\&=
	\dsum_{n=1}^{\mug}\ff{\wt{O}_{n}^{(2)}}{(2n+1)^{2}}t^{2n+1}.
\end{align*}
Clearly, $t=1$ gives the sum for (\ref{47}).
Therefore, 
\[
W\k{\ff{1}{6}(\sin^{-1}t)^{3}}
\Biggr|_{t=1}
=
\di\int_{0}^{1}{
\ff{1}{6}(\sin^{-1}u)^{3}
}\,\ff{du}{\rt{1-u^{2}}}
=
\ts{\ff{1}{24}(\sin^{-1} u)^{4}}{1}{0}=
\ff{\p^{4}}{384}.
\]
Similarly, we have 
\[
W\k{\ff{2}{3}(\sin^{-1}t)^{4}}
=
\ff{\p}{2} \dsum_{n=1}^{\mug}\ff{H_{n-1}^{(2)}}{n^{2}}t^{2n}
\]
so that 
\begin{align*}
	W\k{\ff{2}{3}(\sin^{-1}t)^{4}}
	\Biggr|_{t=1}
	&=
\di\int_{0}^{1}{
\ff{2}{3}(\sin^{-1}u)^{4}
}\,\ff{du}{\rt{1-u^{2}}}
=
\ts{\ff{2}{15}(\sin^{-1} u)^{5}}{1}{0}=
\ff{\p^{5}}{240}.
\end{align*}
We conclude that 
\[
\dsum_{n=1}^{\mug}\ff{H_{n-1}^{(2)}}{n^{2}}
=
\ff{2}{\p} \k{\ff{\p^{5}}{240}}=\ff{\p^{4}}{120}.
\]
\end{proof}
\begin{rmk}\hf
\begin{enumerate}
	\item 
(\ref{47}) is a variation of De Doelder's formula 
$\tsum_{n=1}^{\mug}\tff{O_{n}^{(2)}}{n^{2}}=\tff{\p^{4}}{32}$ 
\cite[p.1196 (13)]{bobo} 
and (\ref{48}) gives another proof of 
$\tsum_{n=1}^{\mug}\tff{H_{n}^{(2)}}{n^{2}}=\tff{7}{4}\z(4)$ 
\te{\cite[p.286]{valean}} 
because 
\[
\dsum_{n=1}^{\mug}\ff{H_{n}^{(2)}}{n^{2}}
=
\dsum_{n=1}^{\mug}
\k{\ff{H_{n-1}^{(2)}}{n^{2}}
+\ff{1}{n^{4}}
}
=
\ff{3}{4}\z(4)+\z(4)
=
\ff{7}{4}\z(4).
\]
	\item After preparation of the manuscript, Christophie Vignat kindly told me that recently Guo-Lim-Qi (2021) \cite{glq} completely described 
	Maclaurin series of integer powers of arcsin.
Also, we found that J.M. Borwein-Chamberland (2007)  \cite{boch} have 
already obtained the same result.
\end{enumerate}
\end{rmk}

\subsection{Integral evaluation}

As byproduct of our discussions, we find evaluation of many integrals with known special values of $\dil(t), \tri(t)$. 
Here, we record several examples. 
Let $\ph=\tff{1+\rt{5}}{2}$ be the golden ratio. Observe that 
\[
\phi^{-1}=\ff{\rt{5}-1}{2}, \q 
\phi^{-2}=\ff{3-\rt{5}}{2}.
\]
We write $\log^{2}x$ for $(\logx)^{2}$.
Note that \[
\log^{2}(\ph^{-1})=
(\log(\ph^{-1}))^{2}=
(-\log(\ph))^{2}=
(\log(\ph))^{2}=
\log^{2}(\ph).
\]

\begin{fact}[\cite{lewin}]
\begin{eqnarray}
\dil(\ph^{-1})=
-\log^{2}(\phi)+\ff{\p^{2}}{10}.
\label{51}
\end{eqnarray}
\begin{eqnarray}
\dil(\ph^{-2})=
-\log^{2}(\phi)+\ff{\p^{2}}{15}.
\label{52}
\end{eqnarray}
\begin{eqnarray}
\tri(\ph^{-2})=
\ff{4}{5}\z(3)-\ff{2\p^{2}}{15}\log\phi+
\ff{2}{3}\log^{3}\phi.
\label{53}
\end{eqnarray}
\begin{eqnarray}
\dil\k{\f12}=\ff{\p^{2}}{12}-\f{\,1\,}{2}\log^{2}2.
\label{54}
\end{eqnarray}
\begin{eqnarray}
\te{Li}_{3}\k{\ff{1}{2}}=
\ff{7}{8}\z(3)-\ff{\p^{2}}{12}\log2+\ff{1}{6}\log^{3}2.
\label{55}
\end{eqnarray}

\end{fact}

\begin{cor}
\begin{eqnarray}
\di\int_{0}^{1}
{
\ff{\sin^{-1}(\phi^{-1}u)}{\rt{1-u^{2}}}
\,{du}}
=
-\ff{3}{4}\log^{2}(\ph)+\ff{\p^{2}}{12}.
\label{56}
\end{eqnarray}
\begin{eqnarray}
\di\int_{0}^{1}{
\ff{\tf12
({\sin^{-1}\tf{u}{\rt{2}}})
^{2}
}{\rt{1-u^{2}}}
}\,du=
\ff{\p}{8}
\k{
\ff{\p^{2}}{12}-\ff12\log^2{2}
}.
\label{57}
\end{eqnarray}
\begin{eqnarray}
\ff{16}{\p}
\di\int_{0}^{1}{
\ff{1}{2}\k{\sin^{-1}\phi^{-1}x}
^{2}
\ff{\cos^{-1}x}{x}
}\,dx=
\ff{4}{5}\z(3)-\ff{2\p^{2}}{15}\log\phi+
\ff{2}{3}\log^{3}\phi.
\label{58}
\end{eqnarray}
\begin{eqnarray}
\di\int_{0}^{1}{
\ff12
\k{{\sinh^{-1}
\phi^{-1/2}u
}}
^{2}
\ff{du}
{\rt{1-u^{2}}}}
\,=
\ff{\p}{2}
\k{
-\ff18\log^2{\ph}+
\ff{\p^{2}}{60}
}.
\label{59}
\end{eqnarray}
\begin{eqnarray}
\ff{16}{\p}
\di\int_{0}^{1}{
\ff{1}{2}\k{\sin^{-1}\ff{x}{\rt{2}}}
^{2}
\ff{\cos^{-1}x}{x}
}\,dx=
\ff{7}{8}\z(3)-\ff{\p^{2}}{12}\log2+\ff{1}{6}\log^{3}2.
\label{60}
\end{eqnarray}
\end{cor}

\begin{proof}
\[
\di\int_{0}^{1}{
\ff{\sin^{-1}
(\phi^{-1}u)}{\rt{1-u^{2}}}
\,{du}}
=
\chi_{2}(\phi^{-1})=\dil(\phi^{-1})-\f{1}{4}\dil(\ph^{-2})
\]
\[
=
\k{-\log^{2}(\ph)+\ff{\pi^{2}}{10}}
-\ff{\,1\,}{4}
\k{-\log^{2}(\ph)+\ff{\pi^{2}}{15}}
=
-\ff{3}{4}\log^{2}(\ph)+\ff{\p^{2}}{12}.	
\]
%
\[
W(\tf12(\sin^{-1}t)^{2})
\bigr|_{t=1/\rt{2}}=
\ff{\p}{8}\dil\k{\f12}
=
\ff{\p}{8}\k{\ff{\p^{2}}{12}-\f{\,1\,}{2}\log^{2}2}.
\]
(\ref{34}) for $t=\phi^{-1}$ with (\ref{53})
 gives (\ref{58}).
\[
W(\tf12(\sinh^{-1}t)^{2})
\bigr|_{t=\phi^{-1/2}}=
\ff{\p}{2}\k{
\ff{1}{4}\dil(\ph^{-1})-
\ff{1}{8}\dil(\ph^{-2})
}
=
\ff{\p}{2}\k{-
\f{\,1\,}{8}\log^{2}\ph+
\ff{\p^{2}}{60}}.
\]
Finally,  (\ref{34}) for $t=1/\rt{2}$ with (\ref{55}) 
gives (\ref{60}).
\end{proof}

\kill{\np}
\section{Concluding remarks}

Here, we record several remarks for our future research.

\begin{enumerate}
	\item For $0\le \al\le 1$, define a generalized Wallis operator 
\[
W_{\al}f(t)=
	\di\int_{0}^{\al}{f(tu)}\,\ff{du}{\rt{1-u^{2}}}
\]
so that we can deal with more general integrals. 
Study $W_{\al}$, particularly for $\al=1/2, \rt{2}/2, \rt{3}/2$.
\item Can we show any inequality for $\tri(t), \chi_3(t)$ and $\te{Ti}_{3}(t)$ in a similar way?
	\item Discuss $(\sinh^{-1}t)^{3}$, $(\sinh^{-1}t)^{4}$ and related Euler sums.
	\item Wolfram alpha \cite{wolfram} says that 	
\[
\di\int_{0}^{1}{\ff{(\sin^{-1}x)^{3}}{x}}\,dx
=
\di\int_{0}^{\p/2}{u^{3}\cot u}\,du
=\ff{\p^{3}}{8}\log2-\ff{9}{16}\pi\z(3),
\]
\[
\di\int_{0}^{1}{\ff{(\sin^{-1}x)^{4}}{x}}\,dx
=
\di\int_{0}^{\p/2}{u^{4}\cot u}\,du
=
\ff{1}{32}\k{-18\p^{2}\z(3)+93\z(5)+2\p^{4}\log2}.
\]
It should be possible to describe such integrals as  certain infinite sums with or without numbers $w_{2n}$.  We plan to study those details in subsequent publication.
\item It is interesting that (\ref{38}) happens to be quite similar to
\[
\di\int_{0}^{1}{\ff{\tan^{-1}x\cot^{-1}x}{x}}\,dx
=\ff{7}{8}\z(3).
\]
Not often this result appears in this form in the literature, though. 
Now, let us see how we evaluate this integral. 
Let 
\[
I=
\di\int_{0}^{1}{\ff{\tan^{-1}x\cot^{-1}x}{x}}\,dx,
\]
\[
I_{1}=
\di\int_{0}^{1}{\ff{\tan^{-1}x}{x}}\,dx,
\q 
I_{2}=
\di\int_{0}^{1}{\ff{(\tan^{-1}x)^{2}}{x}}\,dx.
\]
Then 
\begin{align*}
	I&=\di\int_{0}^{1}{\ff{\tan^{-1}x\cot^{-1}x}{x}}\,dx
	\\&=\di\int_{0}^{1}{\ff{\tan^{-1}x\k{
\f{\p}{2}-\tan^{-1}x}
}{x}}\,dx
=
{\ff{\p}{2}I_{1}-I_{2}}.
\end{align*}
We can compute $I_{1}$ and $I_{2}$ as follows.
\begin{align*}
	I_{1}&=\di\int_{0}^{1}{\ff{\tan^{-1}x}{x}}\,dx
	=\di\int_{0}^{1}{
\dsum_{n=0}^{\mug}\ff{(-1)^{n}}{2n+1}x^{2n}
}\,dx
	\\&=\dsum_{n=0}^{\mug}\ff{(-1)^{n}}{2n+1}
\di\int_{0}^{1}{x^{2n}}\,dx=
\dsum_{n=0}^{\mug}\ff{(-1)^{n}}{(2n+1)^{2}}=
G.
\end{align*}
For $I_{2}$, recall from Fourier analysis that 
\[
\log\k{\tan\ff{y}{2}}=-2
\dsum_{n=0}^{\mug}\ff{1}{2n+1}\cos(2n+1)y, \q 0<y<\pi.
\]
It follows that 
\begin{align*}
	I_{2}&=\di\int_{0}^{1}{\ff{(\tan^{-1}x)^{2}}{x}}\,dx
\xeq[y=2\atx]{}
\ff{\,1\,}{4}
\di\int_{0}^{\p/2}{\ff{y^{2}}{\sin y}}\,dy
	\\&=\ff{\,1\,}{4}
\k{\ts{y^{2}
\log\k{\tan \ff{y}{2}}
}{\p/2}{0}
-
\di\int_{0}^{\p/2}{
2y\log\k{\tan \ff{y}{2}}
}\,dy
}
	\\&=-\ff{1}{2}
\di\int_{0}^{\p/2}{
y
\k{-2
\dsum_{n=0}^{\mug}\ff{1}{2n+1}
\cos(2n+1)y
}
}\,dy
	\\&=\dsum_{n=0}^{\mug}
\ff{1}{2n+1}
\di\int_{0}^{\p/2}{y\cos(2n+1)y}\,dy
	\\&=\dsum_{n=0}^{\mug}
\ff{1}{2n+1}
\k{\ts{
y\ff{\sin(2n+1)y}{2n+1}
}{\p/2}{0}-
\di\int_{0}^{\p/2}{
\ff{\sin(2n+1)y}{2n+1}
}\,dy
}
	\\&=\dsum_{n=0}^{\mug}
\k{\ff{\,\pi\,}{2}\ff{(-1)^{n}}{(2n+1)^{2}}
-\ff{1}{(2n+1)^{3}}}
=\ff{\pi G}{2}-\ff{7}{8}\z(3).
\end{align*}
Finally, we see 
\[
I=
{
\ff{\p G}{2}-
\k{\ff{\pi G}{2}-\ff{7}{8}\z(3)
}
}=\ff{7}{8}\z(3).
\]
Open question:What if we replace $\tan^{-1}$ by $\tanh^{-1}$?
\end{enumerate}

In this article, we encountered many integral representations for dilogarithm, trilogarithm and hence $\z(2)$, the Catalan constant $G$ and $\z(3)$ as a reformulation of Boo Rim Choe (1987) \cite{boo}, Ewell (1990) \cite{ewell} and Williams-Yue (1993) \cite{yue} on the inverse sine function. 
As an application, we also proved new Euler sums. Indeed, there are subsequent results on multiple zeta and $t$-values $\z(3, 2, \cd, 2)$, $t(3, 2, \cd, 2)$ as Hoffman and Zagier discussed in \cite{hoffman, zagier2}. 
We will write them with more details at another  opportunity.
\begin{center}
\textbf{Acknowledgment.}
\end{center}

\begin{center}
This research arose from Iitaka online seminar in 2020-2021. 
The author would like to thank the organizer 
Shigeru Iitaka, Shunji Sasaki and Kouichi Nakagawa for fruitful discussions. Satomi Abe as well as 
Michihito Tobe sincerely supported him for his writing. 
He also thanks the anonymous reviewer for helpful comments to improve the manuscript.
\end{center}

\end{document}